\documentclass[12pt,leqno]{article}
\usepackage[pdftex]{graphicx}
\usepackage{amsfonts,amssymb,amsmath}
\usepackage{ccaption,mathrsfs,amsmath,amssymb,amsthm}
\usepackage{amsmath, amscd, amsfonts,  tikz, times}

\newcounter{conjecture}
\newcounter{remark}

\newcommand{\eqnsection}{
    \renewcommand{\theequation}{\thesection.\arabic{equation}}
    \makeatletter
    \csname @addtoreset\endcsname{equation}{section}
    \makeatother}
\newtheorem{theorem}{Theorem}

\newtheorem{lemma}{Lemma}

\newcommand{\dd}{\delta}

\newcommand{\lar}{\longrightarrow}
\newcommand{\eps}{\varepsilon}

\newcommand{\CC}{\mathbb{C}}

\def \be{\begin{equation}}
\def \ee{\end{equation}}
\def \bt{\begin{theorem}}
\def \et{\end{theorem}}
\def \bea{\begin{eqnarray}}
\def \eea{\end{eqnarray}}
\def \bas{\begin{eqnarray*}}
\def \eas{\end{eqnarray*}}

\newcommand {\rrr}[1]{(\ref{#1})}

\def \sech{\mbox{sech}}
\def \csch{\mbox{csch}}
\def \aa{\alpha}
\def \bb{\beta}

\def \ga{\gamma}

\def \th{\theta}

\def \ff{\infty}

\def \DD{{\mathbb D}}

\def \HH{{\mathbb H}}

\def \RR{{\mathbb R}}

\def \CCC{{\cal C}}

\def \({\left(}
\def \){\right)}

\def \vski{\vspace{12pt}}

\def \bc{\begin{center} }
\def \ec{\end{center} }
\def \bs{\begin{slide} }
\def \es{\end{slide} }

\def\square{{\vcenter{\vbox{\hrule height.3pt
         \hbox{\vrule width.3pt height5pt \kern5pt
            \vrule width.3pt}
         \hrule height.3pt}}}}
\def\qed{{\hfill $\Box$ \bigskip}}

\newcounter{ccases}
\newcommand{\ccases}[1]{\begingroup \refstepcounter{ccases} {\bf \fontsize{12}{12}\selectfont Example \theccases }  \label{#1}\endgroup}

\eqnsection
\begin{document}

\title{On the distribution of planar Brownian motion at stopping times.}

\author{
\begin{tabular}{c}
\textit{Greg Markowsky} \\
Monash University \\
Victoria 3800, Australia \\
gmarkowsky@gmail.com
\end{tabular}}

\bibliographystyle{amsplain}

\maketitle \eqnsection \setlength{\unitlength}{2mm}

\begin{abstract}

A simple extension is given of the well-known conformal invariance of harmonic measure in the plane. This equivalence depends on the interpretation of harmonic measure as an exit distribution of planar Brownian motion, and extends conformal invariance to analytic functions which are not injective, as well as allowing for stopping times more general than exit times. This generalization allow considerations of homotopy and reflection to be applied in order to compute new expressions for exit distributions of various domains, as well as the distribution of Brownian motion at certain other stopping times. An application of these methods is the derivation of a number of infinite sum identities, including the Leibniz formula for $\pi$ and the values of the Riemann $\zeta$ function at even integers.

\vski


\vski

\end{abstract}

\section{Introduction and primary methods}

It is well known that harmonic measure on domains in $\CC$ can be interpreted in terms of exit distributions of planar Brownian motion. The conformal invariance of harmonic measure therefore implies a conformal invariance principle for exit distributions as well, and this principle can also be deduced directly from L\'evy's theorem on the conformal invariance of Brownian motion. This invariance allows in many instances for exit distributions to be calculated on simply connected domains for which a conformal equivalence with the disk is known. However, L\'evy's theorem in fact does not require maps to be injective, permitting general nonconstant analytic functions as well. We will show how the conformal invariance of the exit distribution of Brownian motion can be extended to nonconstant analytic functions, as well as to more general stopping times. This allows us to derive new expressions for many exit times and to calculate the distributions of Brownian motion at the exit time of certain non-simply connected domain, as well as at certain stopping times which are not exit times. We will give a number of illustrative examples, and show how certain identities can result from the appropriate choice of stopping times. Perhaps most notably, we will see a number of different ways in which the values for $\sum_{n=1}^{\ff} \frac{1}{n^{2m}}$ and $\sum_{n=1}^{\ff} \frac{(-1)^{n-1}}{n^{2m+1}}$ can be deduced.

\vski

L\'evy's theorem is as follows (see \cite{durBM} or \cite{revyor} for a proof).

\begin{theorem} \label{holinv}
Let $f$ be analytic and nonconstant on a domain $U$, and let $a \in U$. Let $B_t$ be a Brownian motion starting at $a$, and $\tau$ a stopping time such that the set of Brownian paths $\{B_t: 0 \leq t \leq \tau\}$ lie within $U$ a.s. Then the process $f(B_{t})$ stopped at $\tau$ is a time-changed Brownian motion.
\end{theorem}

It should be noted that the time change referenced in the theorem can be expressed explicitly, but is not important for our purposes. Let $\ga$ be a smooth curve parameterized by arclength, $B_t$ a Brownian motion starting at $a$, and $\tau$ a stopping time such that $B_\tau \in \ga$ a.s. $\rho^a_\tau(w) ds$ will denote the density of $B_\tau$ on $\ga$, when it exists, with $ds$ denoting the arclength element. We then have the following identity, valid for any measurable subset $A$ of $\ga$:

\begin{equation} \label{}
P_a(B_\tau \in A) = \int_{A} \rho^a_\tau(\ga_s) ds.
\end{equation}

Now, Levy's Theorem provides us with a Brownian motion $\hat B_t$ which is a time-change of $f(B_t)$ and a stopping time $\hat \tau$, which is the image under the time change of $\tau$, so that $\hat B_{\hat \tau} \in f(\ga)$ a.s. We will also use the notation $\rho^a_{\hat \tau}(w) ds$ to denote the density of $\hat B_{\hat \tau} = f(B_\tau)$. Our method of projection is contained in the following theorem.

\begin{theorem} \label{massage}
Let $U$ be a domain, and suppose $f$ is a function analytic on $U$. Let $B_t$ be a Brownian motion starting at $a$, and $\tau$ a stopping time such that the set of Brownian paths $\{B_t: 0 \leq t \leq \tau\}$ lie within $U$ a.s. Suppose that $\ga$ is a smooth curve in $U$ such that $B_\tau \in \ga$ a.s. Then for any $a \in V$ and $w \in f(\ga)$ we have

\begin{equation} \label{cherry}
\rho^{f(a)}_{\hat \tau} (w)ds = \sum_{z \in f^{-1}(w)\cap \ga} \frac{\rho^a_\tau (z)}{|f'(z)|} ds.
\end{equation}

\end{theorem}

The proof is essentially immediate, since any Brownian path which finishes at $f^{-1}(w) \cap \ga$ at time $\tau$ will be mapped under $f$ to a path finishing at $w$. The $|f'(z)|$ in the denominator on the right side of \rrr{cherry} is the scaling factor required for the change in the arclength element mapped under the analytic function $f$. We also remark that in the case that $f$ is conformal, it is often easier to use $g = f^{-1}$, and \rrr{cherry} becomes

\begin{equation} \label{cherry2}
\rho^{a}_{\hat \tau} (w)ds = \rho^{g(a)}_{\tau} (g(w))|g'(w)| ds.
\end{equation}

In the next section we will proceed through a series of examples which illustrate the use of the theorem.

\section{Examples}

The first two examples, concerning the disk and half-plane, are certainly known but are included for completeness and for their use in the later examples. In what follows, we will use the notation $T_U$ to denote the exit time of any domain $U$; that is, $T_U = \inf\{t \geq 0 | B_t \in U^c\}$.

\vski

\ccases{disk}(disk) There is only one exit distribution in $\CC$ which is obvious: that of a disk in which the Brownian motion starts at the center. With $\DD = \{|z|<1\}$, rotational invariance shows immediately that $\rho^0_{T_\DD}(e^{i \th})ds = \frac{ds}{2\pi}$. For $a \in \DD$ consider the M\"obius transformation

\begin{equation} \label{diskmap}
\psi_a(z) = \frac{z-a}{1-\bar a z}.
\end{equation}

It is well-known that $\psi_a$ is a conformal self-map of $\DD$ sending $a$ to 0. Using $\psi'_a = \frac{1-|a|^2}{(1-\bar a z)^2}$ and \rrr{cherry2} we obtain the following identity:

\begin{equation} \label{mae}
\rho^a_{T_\DD}(e^{i \th})ds = \frac{1}{2\pi} \frac{1-|a|^2}{|1-\bar a e^{i\th}|^2} ds.
\end{equation}

We can also calculate the exit density from disks with radius other than one, as well as the hitting density of the circle when the Brownian motion begins at a point outside the circle. Let $m\DD = \{|z|<m\}$ and $m\DD^c = \{|z|>m\}$; Note that $T_{m\DD}$ and $T_{m\DD^c}$ both signify the first hitting time of $\{|z|=m\}$, but the former is of interest when the initial point of the Brownian motion has modulus less than $m$, and the latter when it is greater. If $|a|<m$, then we can project \rrr{mae} using the map $z \lar mz$ to obtain

\begin{equation} \label{}
\rho^a_{T_{m\DD}}(me^{i \th})ds = \frac{1}{m}\rho^{a/m}_{T_{\DD}}(e^{i \th})ds = \frac{1}{2\pi m} \frac{1-|\frac{a}{m}|^2}{|1-\frac{\bar a e^{i \th}}{m}|^2}ds = \frac{1}{2\pi m} \frac{m^2-|a|^2}{|m-\bar a e^{i \th}|^2}ds,
\end{equation}

while if $|a|>m$ then we use the map $z \lar \frac{m}{z}$ to get

\begin{equation} \label{}
\rho^a_{T_{m\DD}}(me^{i \th})ds = \frac{1}{m}\rho^{m/a}_{T_{\DD}}(e^{-i \th})ds = \frac{1}{2\pi m} \frac{1-|\frac{m}{a}|^2}{|1-\frac{m}{\bar a} e^{-i \th}|^2}ds = \frac{1}{2\pi m} \frac{|a|^2-m^2}{|a-m e^{i \th}|^2}ds.
\end{equation}

As a side note, the Poisson Integral Formula for harmonic functions in the disk $\DD$ can be derived from \rrr{mae}, as we have by Dynkin's formula for harmonic $h$ (see \cite{fima})

\begin{equation} \label{}
h(a) = E_a[h(B_{T_\DD})] = \frac{1}{2\pi} \int_{0}^{2\pi} h(e^{i\th})\frac{1-|a|^2}{|1-\bar a e^{i\th}|^2} d\th.
\end{equation}

The reader may check that setting $a = re^{it}$ and performing a few simple manipulations yield a more standard form of the formula.

\vski

\ccases{halfplane} (half-plane) We can also easily calculate the exit distribution of a half-plane, as follows. Let $T_\HH$ be this exit time. The conformal map taking $\HH = \{y>0\}$ to $\DD$ is given by $f(z) = \frac{z-i}{z+i}$, with $f'(z) = \frac{2i}{(z+i)^2}$. We obtain

\begin{equation} \label{}
\rho^i_{T_\HH}(x)ds = \frac{1}{\pi} \frac{1}{1+x^2}ds.
\end{equation}

To find the distribution from a more general point $a = u+vi$ use the map $f(z) = u+vz$, which fixes $\HH$ and maps $i$ to $a$, to obtain

\begin{equation} \label{}
\rho^a_{T_\HH}(x)ds = \frac{1}{\pi} \frac{v}{v^2+(x-u)^2}ds.
\end{equation}

As in the case of the unit disk, this distribution leads via Dynkin's formula to the Poisson Integral Formula for the upper half-plane.

\vski

\ccases{punctdisk} (punctured disk) Now let $V = \{0<|z|<1\}$ be the punctured disk. We can calculate the exit distribution of $V$ by projecting from $\HH$ via the covering map $f(z) = e^{iz}$ (note that $|f'(z)| = 1$ on $\dd \HH$). For any $a \in V$ we obtain

\begin{equation} \label{page}
\rho^a_{T_V}(e^{i \th})ds= \sum_{k=-\ff}^{\ff} \frac{-\ln |a|}{\pi ((\ln |a|)^2 + (arg(a)-(\th+2\pi k))^2)}ds.
\end{equation}

Planar Brownian motion does not see points, i.e. $P_a(B_t = 0 \mbox{ for some }t \geq 0) = 0$. Thus, the exit distribution for the disk and the punctured disk agree. The difference then between the expressions \rrr{mae} and \rrr{page} is that each term in the sum in \rrr{page} corresponds to a different homotopy class of paths in the punctured disk terminating at $e^{i\th}$, while \rrr{mae} does not differentiate between the homotopy classes. Furthermore, equating \rrr{mae} and \rrr{page}, assuming $a \in (0,1)$ for simplicity, gives the identity

\begin{equation} \label{effyou}
\sum_{k=-\ff}^{\ff} \frac{-\ln a}{\pi ((\ln a)^2 + (\th+2\pi k)^2)} = \frac{1}{2\pi} \frac{1-a^2}{|1-a e^{i\th}|^2} = \frac{1-a^2}{2\pi(1+a^2-2a\cos \th)}.
\end{equation}

This identity can be manipulated into a more easily recognized identity, as follows. Divide both sides by $-\ln a$ and simplify. This gives

\begin{equation} \label{cthg}
\sum_{k=-\ff}^{\ff} \frac{1}{\pi((\ln a)^2 + (\th + 2\pi k)^2) } = \frac{1-a^2}{2(-\ln a)(1+a^2-2a\cos \th)}.
\end{equation}

Assuming $\th \neq 0$, we can now let $a \nearrow 1$ using $\lim_{a \lar 1} \frac{1-a^2}{\ln a} = -2$, and obtain

\begin{equation} \label{cthg2}
\sum_{k=-\ff}^{\ff} \frac{1}{(\th + 2\pi k)^2 } = \frac{1}{2(1-\cos \th)}.
\end{equation}

Subtract $\frac{1}{\th^2}$ from both sides (the term $k=0$), take the limit as $\th \lar 0$ using $\lim_{\th \lar 0} \frac{1}{2(1-\cos \th)} - \frac{1}{\th^2} = \frac{1}{12}$, and simplify. We obtain

\begin{equation} \label{jill}
\sum_{k=1}^{\ff} \frac{1}{k^2} = \frac{\pi^2}{6}.
\end{equation}

This is Euler's celebrated Basel sum, which has many other existing proofs, including a different probabilistic proof making use of planar Brownian motion (see \cite{meecp}). Note also that if we differentiate \rrr{cthg2} $2m-2$ times and let $\th \lar 0$ we will be able to obtain the well-known values of $\sum_{k=1}^{\ff} \frac{1}{k^{2m}}$. Furthermore, if we take $a = e^{-1}$, then $\th = 0, \pi$ successively, we obtain two identities, which may be added to obtain a third as follows.

\begin{equation} \label{mapleton}
\begin{gathered}
\sum_{k=-\ff}^{\ff} \frac{1}{\pi (1 + (2\pi k)^2)} = \frac{1}{2}\coth(\frac{1}{2}), \\
\sum_{k=-\ff}^{\ff} \frac{1}{\pi (1 + (\pi+2\pi k)^2)} = \frac{1}{2}\tanh(\frac{1}{2}), \\
\sum_{k=-\ff}^{\ff} \frac{1}{\pi (1 + (\pi k)^2)} = \coth(1).
\end{gathered}
\end{equation}

The final identity in \rrr{mapleton} is a standard identity which arises as an example of several different techniques for summing series, for instance using the residue theorem \cite[Ch. 7]{schaum}.

\vski

\ccases{strip} (infinite strip) We now calculate the exit distribution on an infinite strip. We let $W = \{-1 < Re(z) < 1\}$ and take the starting point of the Brownian motion $a$ to lie on the real interval $(-1,1)$; it is clear that the distribution with any other starting point can be obtained from this merely by translation. Apply our theorem to the function $\tan (\frac{\pi}{4} z)$, which maps $W$ conformally to $\DD$, to get

\begin{equation} \label{is1}
\rho^a_{T_W}(\pm 1+yi)ds = |\frac{\pi}{4}\sec^2(\frac{\pi}{4}(\pm 1+yi)))|\rho^{\tan (\frac{\pi}{4} a)}_{T_\DD}(\tan (\pm \frac{\pi}{4} + \frac{\pi}{4}yi))ds.
\end{equation}

We can simplify

\begin{equation} \label{}
\begin{split}
|\frac{\pi}{4}\sec^2(\frac{\pi}{4}(\pm 1+yi)))| & = \pi |e^{-\frac{\pi}{4}y+ \frac{\pi}{4}i} + e^{\frac{\pi}{4}y- \frac{\pi}{4}i}|^{-2} \\
& = \pi |\frac{1}{\sqrt{2}}(e^{-\frac{\pi}{4}y}+ e^{\frac{\pi}{4}y}) + \frac{i}{\sqrt{2}}(e^{-\frac{\pi}{4}y}- e^{\frac{\pi}{4}y})|^{-2} \\
& = \frac{\pi}{2}(\cosh^2 \frac{\pi}{4}y+\sinh^2 \frac{\pi}{4}y)^{-1} \\
& = \frac{\pi}{2}(\cosh \frac{\pi}{2}y)^{-1}.
\end{split}
\end{equation}

Furthermore, using the distribution calculated in Example \ref{disk}, we have

\begin{equation} \label{}
\rho^{\tan (\frac{\pi}{4} a)}_{T_\DD}(\tan (\pm \frac{\pi}{4} + \frac{\pi}{4}yi))ds = \frac{1-\tan^2(\frac{\pi}{4} a)}{2\pi|1-\tan (\frac{\pi}{4} a)\tan (\pm \frac{\pi}{4} + \frac{\pi}{4}yi)|^2}ds.
\end{equation}

Using $\tan(\aa + \bb) = \frac{\tan \aa + \tan \bb}{1-\tan \aa \tan \bb}$ and the identity $\tan(\frac{\pi}{4}yi)= i \tanh(\frac{\pi}{4}y)$ gives

\begin{equation} \label{}
\begin{split}
\tan (\pm \frac{\pi}{4} + \frac{\pi}{4}yi) &= \frac{\pm 1 + i \tanh(\frac{\pi}{4}y)}{1 \mp i \tanh(\frac{\pi}{4}y)} = \frac{\pm(1 \pm i \tanh(\frac{\pi}{4}y))^2}{1 + \tanh^2(\frac{\pi}{4}y)} \\
& = \frac{\pm(1 - \tanh^2(\frac{\pi}{4}y)) + 2i\tanh(\frac{\pi}{4}y) }{1 + \tanh^2(\frac{\pi}{4}y)},
\end{split}
\end{equation}

which yields

\begin{equation} \label{is2}
\rho^{\tan (\frac{\pi}{4} a)}_{T_\DD}(\tan (\pm \frac{\pi}{4} + \frac{\pi}{4}yi))ds = \frac{1-\tan^2(\frac{\pi}{4} a)}{2\pi((1\mp\tan (\frac{\pi}{4} a)\frac{1 - \tanh^2(\frac{\pi}{4}y)}{1 + \tanh^2(\frac{\pi}{4}y)})^2+(\tan (\frac{\pi}{4} a)\frac{2\tanh(\frac{\pi}{4}y) }{1 + \tanh^2(\frac{\pi}{4}y)})^2)}ds.
\end{equation}

Combining \rrr{is1}-\rrr{is2} gives

\begin{equation} \label{is3}
\rho^a_{T_W}(\pm 1+yi)ds =  \frac{\sech (\frac{\pi}{2}y)(1-\tan^2(\frac{\pi}{4} a))}{4((1\mp\tan (\frac{\pi}{4} a)\frac{1 - \tanh^2(\frac{\pi}{4}y)}{1 + \tanh^2(\frac{\pi}{4}y)})^2+(\tan (\frac{\pi}{4} a)\frac{2\tanh(\frac{\pi}{4}y) }{1 + \tanh^2(\frac{\pi}{4}y)})^2)}ds.
\end{equation}

There is another method for calculating $\rho^a_{T_W}(\pm 1+yi)ds$ which uses a form of the reflection principle. This method will be applied in other examples below, namely Examples \ref{halfstrip} and \ref{rectangle}, and we prove its validity carefully for this example, while in the later ones merely referencing this one. Let us define $\tau(b) = \inf\{t \geq 0: Re(B_t) = b\}$. It is clear with comparison with the half-plane example above that $\rho^a_{\tau(b)}(b+yi)ds = \frac{1}{\pi}\frac{|a-b|}{|a-b|^2+y^2}ds$. It will turn out that

\begin{equation} \label{tats}
\begin{gathered}
\rho^a_{T_W}(1+yi)ds = \rho^a_{\tau(1)}(1+yi)ds - \rho^a_{\tau(-3)}(-3+yi)ds + \rho^a_{\tau(5)}(5+yi)ds - \rho^a_{\tau(-7)}(-7+yi)ds + \ldots, \\
\rho^a_{T_W}(-1+yi)ds = \rho^a_{\tau(-1)}(1+yi)ds - \rho^a_{\tau(3)}(3+yi)ds + \rho^a_{\tau(-5)}(-5+yi)ds - \rho^a_{\tau(7)}(7+yi)ds + \ldots.
\end{gathered}
\end{equation}

We will prove the first equation in \rrr{tats}, and for the proof of this it will help to isolate several lemmas. Let us extend the definition of $\tau$ by recursively defining $\tau(b_1, b_2, \ldots, b_{n+1}) = \inf\{t \geq \tau(b_1, \ldots, b_{n}): Re(B_t) = b_{n+1}\}$; that is, $\tau(b_1, b_2, \ldots, b_{n})$ is the first time at which $Re(B_t)$ has visited the sequence $b_1, b_2, \ldots , b_n$ in order. The first lemma should be clear upon consideration, and we therefore omit the proof.

\begin{lemma} \label{still}
On the event $\{\tau(b_1) < \tau(b_2)\}$, we have $\tau(b_1,b_2, \dots, b_n) = \tau(b_2, \dots, b_n)$.
\end{lemma}

That is, on the set of all Brownian paths which hit $b_1$ before $b_2$, we can drop $b_1$ from the start of the sequence without changing the value of the stopping time. The following is an immediate consequence.

\begin{lemma} \label{mang}
For $A \subseteq \RR$, we have

\begin{equation} \label{alive}
\begin{split}
&\Big( 1_{\{Im(B_{\tau(1)}) \in A\}} - 1_{\{Im(B_{\tau(-1,1)}) \in A\}}\Big) + \Big(1_{\{Im(B_{\tau(1,-1,1)}) \in A\}} - 1_{\{Im(B_{\tau(-1,1,-1,1)}) \in A\}}\Big) \\
&+ \ldots + \Big(1_{\{Im(B_{\tau(1,-1, \ldots, -1, 1)}) \in A\}} - 1_{\{Im(B_{\tau(-1,1,-1,\ldots, 1, 1)}) \in A\}}\Big) \\
& \qquad \leq 1_{\{Re(B_{T_W})=1,Im(B_{T_W}) \in A\}} \\
& \qquad \leq 1_{\{Im(B_{\tau(1)}) \in A\}} -\Big( 1_{\{Im(B_{\tau(-1,1)}) \in A\}} - 1_{\{Im(B_{\tau(1,-1,1)}) \in A\}}\Big) \\
& \qquad  \qquad - \ldots - \Big( 1_{\{Im(B_{\tau(-1,1, \ldots ,-1,1)}) \in A\}} - 1_{\{Im(B_{\tau(1,-1,1, \ldots ,-1,1)}) \in A\}}\Big).
\end{split}
\end{equation}

\end{lemma}

{\bf Remark:} In the sums in \rrr{alive}, it should be understood that the sequences defining the $\tau$'s alternate and increase in length by one with each successive term.

\vski

{\bf Proof:} We begin by noting that $\min(1_{\{Im(B_{\tau(1)}) \in A\}}, 1_{\{\tau(1) < \tau(-1)\}}) = 1_{\{Re(B_{T_W})=1,Im(B_{T_W}) \in A\}}$. By Lemma \ref{still}, on the set $\{\tau(-1) < \tau(1)\}$ we have $\tau(-1,1,-1,\ldots, -1, 1) = \tau(1,-1, \ldots, -1, 1)$, and thus each positive term in the leftmost sum in \rrr{alive} is canceled by the subsequent negative term and therefore the sum is 0. On the other hand, on the set $\{\tau(1) < \tau(-1)\}$ each negative term in the leftmost sum except the last is canceled by the subsequent positive one, again by Lemma \ref{still}. Thus, on $\{\tau(1) < \tau(-1)\}$, the leftmost side of \rrr{alive} is equal to $1_{\{Im(B_{\tau(1)}) \in A\}} - 1_{\{Im(B_{\tau(-1,1,-1,\ldots, 1, 1)}) \in A\}} \leq 1_{\{Im(B_{\tau(1)}) \in A\}} = 1_{\{Re(B_{T_W})=1,Im(B_{\tau(1)}) \in A\}}$. It follows that the leftmost side is less than or equal to $\min(1_{\{Im(B_{\tau(1)}) \in A\}}, 1_{\{\tau(1) < \tau(-1)\}}) = 1_{\{Re(B_{T_W})=1,Im(B_{T_W}) \in A\}}$. The second inequality follows similarly from Lemma \ref{still}, for on the set $\{\tau(1) < \tau(-1)\}$ the difference inside each set of parentheses on the rightmost side is zero yielding a value of $1_{\{Im(B_{\tau(1)}) \in A\}} \geq 1_{\{Re(B_{T_W})=1,Im(B_{T_W}) \in A\}}$, while on $\{\tau(-1) < \tau(1)\}$ we have $1_{\{Re(B_{T_W})=1,Im(B_{T_W}) \in A\}}=0$, while the right side is equal to $1_{\{Im(B_{\tau(1,-1,1, \ldots ,-1,1)}) \in A\}} \geq 0$. \qed

If we let $A$ be a small interval on the line $\{Re(z)=1\}$ centered at $1+yi$, divide by the length of the interval, and then let this length go to 0 we obtain

\begin{equation} \label{denim}
\rho^a_{T_W}(1+yi)ds = \rho^a_{\tau(1)}(1+yi)ds - \rho^a_{\tau(-1,1)}(1+yi)ds + \rho^a_{\tau(1,-1,1)}(1+yi)ds - \rho^a_{\tau(-1,1,-1,1)}(1+yi)ds + \ldots.
\end{equation}

We should mention that, intuitively, \rrr{denim} is very simple: in order to calculate $\rho^a_{T_W}(1+yi)ds$ we want to count paths which leave $\{Re(z)<1\}$ at $1+yi$, however we want to remove the contribution from paths which strike $\{Re(z)=-1\}$ first, so we consider $\rho^a_{\tau(1)}(1+yi)ds - \rho^a_{\tau(-1,1)}(1+yi)ds$; however, by subtracting $\rho^a_{\tau(-1,1)}(1+yi)ds$ we have subtracted too much, as we have incorrectly subtracted the contribution from paths which hit $\{Re(z)=1\}$ before $\{Re(z)=-1\}$, so we must add $\rho^a_{\tau(1,-1,1)}(1+yi)ds$ to compensate; however, by an analogous argument we have overcompensated, and must therefore subtract $\rho^a_{\tau(-1,1,-1,1)}(1+yi)ds$, and so forth. It remains only to understand how to calculate the density $\rho^a_{\tau(b_1, \ldots, b_{n})}(b_n+yi)ds$.

\begin{lemma} \label{}
For any sequence of real numbers $a=b_0, b_1, \ldots, b_n$ we have

\begin{equation} \label{}
\begin{split}
\rho^a_{\tau(b_1, \ldots, b_n)}(b_n + yi)ds &= \rho^a_{\tau(a + \sum_{j=1}^{n}|b_j-b_{j-1}|)}(a + \sum_{j=1}^{n}|b_j-b_{j-1}| + yi)ds \\
&= \rho^a_{\tau(a - \sum_{j=1}^{n}|b_j-b_{j-1}|)}(a - \sum_{j=1}^{n}|b_j-b_{j-1}| + yi)ds.
\end{split}
\end{equation}
\end{lemma}

{\bf Proof:} By induction on $n$. The case $n=1$ follows from the symmetry of Brownian motion over the line $\{Re(z)=a\}$. Suppose that the result holds for $n$, and consider a sequence $b_1, \ldots, b_n, b_{n+1}$. If $b_n$ lies between $b_{n-1}$ and $b_{n+1}$, then $\tau(b_1, \ldots, b_{n-1}, b_n, b_{n+1}) = \tau(b_1, \ldots, b_{n-1}, b_{n+1})$, since the real part of the Brownian motion must hit $b_n$ in passing from $b_{n-1}$ to $b_{n+1}$, so the result follows from the induction hypothesis (since then also $|b_{n+1}-b_n| + |b_n-b_{n-1}| = |b_{n+1}-b_{n-1}|$). If, on the other hand, $b_n$ does not lie between $b_{n-1}$ and $b_{n+1}$, then it must lie between $b_{n-1}$ and $b_n -(b_{n+1}-b_n)$. However, we must have $\rho^a_{\tau(b_1, \ldots, b_n, b_{n+1})}(b_{n+1} + yi)ds = \rho^a_{\tau(b_1, \ldots, b_n, b_n-(b_{n+1}-b_n))}(b_n-(b_{n+1}-b_n) + yi)ds$, since the reflection principle implies that the process

\begin{equation} \label{}
\tilde B_t = \left \{ \begin{array}{ll}
B_t & \qquad  \mbox{if } t \leq \tau(b_1, \ldots, b_n) \\
b_n - (Re(B_t) - b_n) + i Im(B_t) & \qquad \mbox{if } t > \tau(b_1, \ldots, b_n)\;,
\end{array} \right.
\end{equation}

is also a Brownian motion; this is the reflection of $B_t$ over the line $\{Re(z) = b_n\}$ for $t > \tau(b_1, \ldots, b_n)$. By the same argument as before, $\tau(b_1, \ldots,b_{n-1}, b_n, b_n-(b_{n+1}-b_n)) = \tau(b_1, \ldots,b_{n-1}, b_n-(b_{n+1}-b_n))$, and furthermore $|b_n-(b_{n+1}-b_n) - b_n| = |b_{n+1}-b_n|$, so the result again follows from the induction hypothesis. \qed

Using this lemma, we see that \rrr{denim} reduces to the first equation in \rrr{tats}, with the other equation in \rrr{tats} following by the symmetric argument. Using the expression immediately preceding \rrr{tats}, we see

\begin{equation} \label{tats2}
\begin{gathered}
\rho^a_{T_W}(1+yi)ds = \frac{ds}{\pi} \Big( \frac{1-a}{(1-a)^2+y^2} - \frac{3+a}{(3+a)^2+y^2} + \frac{5-a}{(5-a)^2+y^2} - \frac{7+a}{(7+a)^2+y^2} + \ldots \Big), \\
\rho^a_{T_W}(-1+yi)ds = \frac{ds}{\pi} \Big( \frac{1+a}{(1+a)^2+y^2} - \frac{3-a}{(3-a)^2+y^2} + \frac{5+a}{(5+a)^2+y^2} - \frac{7-a}{(7-a)^2+y^2} + \ldots \Big).
\end{gathered}
\end{equation}

Equating the expression for $\rho^a_{T_W}(1+yi)ds$ in \rrr{tats2} with that in \rrr{is3} gives the identity

\begin{equation} \label{nosering}
\begin{split}
\frac{1}{\pi}\sum_{j=1}^{\ff} \frac{(-1)^{j+1}((2j-1) + (-1)^j a)}{((2j-1) + (-1)^j a)^2 +y^2} = \frac{\sech (\frac{\pi}{2}y)(1-\tan^2(\frac{\pi}{4} a))}{4((1-\tan (\frac{\pi}{4} a)\frac{1 - \tanh^2(\frac{\pi}{4}y)}{1 + \tanh^2(\frac{\pi}{4}y)})^2+(\tan (\frac{\pi}{4} a)\frac{2\tanh(\frac{\pi}{4}y) }{1 + \tanh^2(\frac{\pi}{4}y)})^2)}.
\end{split}
\end{equation}

Needless to say, our lives are considerably simplified by setting $y=0$ or $a=0$. For $a=0$, we obtain

\begin{equation} \label{}
\sum_{j=1}^{\ff} \frac{(-1)^{j+1}(2j-1)}{(2j-1)^2 +y^2} = \frac{\pi}{4}\sech(\frac{\pi}{2}y),
\end{equation}

which can be obtained by other methods, for instance the residue theorem (\cite{schaum}). Taking $y=0$ gives us Leibniz's representation for $\pi$:

\begin{equation} \label{propasi}
\frac{\pi}{4} = 1-\frac{1}{3} + \frac{1}{5} - \frac{1}{7} + \ldots
\end{equation}

Returning to \rrr{tats2}, set now $y=0$ to obtain

\begin{equation} \label{mei}
\sum_{j=1}^{\ff} \frac{(-1)^{j+1}}{(2j-1) + (-1)^j a } = \frac{1}{1-a} - \frac{1}{3+a} + \frac{1}{5-a} - \frac{1}{7-a} + \ldots=  \frac{\pi}{4} \Big(\frac{1+\tan (\frac{\pi}{4} a)}{1-\tan (\frac{\pi}{4} a)}\Big).
\end{equation}

This identity is somewhat unusual, and may be new. It can be manipulated to obtain a number of other identities, as follows. Let $g(a) = \frac{\pi}{4} \Big(\frac{1+\tan (\frac{\pi}{4} a)}{1-\tan (\frac{\pi}{4} a)}\Big)$. Then, if $r$ is a positive integer, differentiating \rrr{mei} $r-1$ times yields the identity

\begin{equation} \label{mei2}
\sum_{j=1}^{\ff} \frac{(r-1)!(-1)^{r(j+1)}}{((2j-1) + (-1)^j a)^r} = g^{(r-1)}(a).
\end{equation}

Note that if $r$ is even then all terms in the sum will be positive, while if $r$ is odd then the sum will be alternating. Setting $a = 0$ in this identity gives the values of all sums of the form $\sum_{j=1}^{\ff} \frac{1}{(2j-1)^{2m}}$ or $\sum_{j=1}^{\ff} \frac{(-1)^j}{(2j-1)^{2m+1}}$, which are well-known with the first few equal to $\sum_{j=1}^{\ff} \frac{(-1)^{j+1}}{(2j-1)} = \frac{\pi}{4}, \sum_{j=1}^{\ff} \frac{1}{(2j-1)^{2}} = \frac{\pi^2}{8}, \sum_{j=1}^{\ff} \frac{(-1)^{j+1}}{(2j-1)^{3}} = \frac{\pi^3}{32}, \sum_{j=1}^{\ff} \frac{1}{(2j-1)^{4}} = \frac{\pi^4}{96}, \ldots$. These values have been known since the time of Euler, and it should be noted that the values of the even-powered identities easily give the values of $\zeta(2k)$, the Riemann $\zeta$ function evaluated at the even integers.

\vski

\rrr{mei2} can also be manipulated into a different identity, as follows. Let $r$ be odd and let us refer to the identity \rrr{mei2} as $I(a)$. Let $q$ be a positive integer, and consider the sum $I(\frac{q-1}{q}) + I(\frac{q-3}{q}) + \ldots + I(\frac{-q+3}{q}) + I(\frac{-q+1}{q})$. The left side of this new identity is a sum, and we will group the terms in this new sum according to their place in the original sums; that is, the first term will be the sum of the first terms in $I(\frac{q-1}{q}), \ldots, I(\frac{-q+1}{q})$, the second term will be the sum of the corresponding second terms (each of which is negative), and so forth. Using $\frac{1}{((2j-1)+ (-1)^j (\frac{p}{q}))^r} = \frac{q^r}{((2j-1)q+ (-1)^j p)^r}$, it may be verified that the first term in the new sum is $q^r$ times the sum of the reciprocals of the $r$-th power of the odd integers from $1$ to $2q-1$, the second term (which is negative) is $q^r$ times the sum of the reciprocals of the $r$-th power of the odd integers from $2q+1$ to $4q-1$, and so forth. We obtain the identity

\begin{equation} \label{armdark}
\begin{split}
\Big( \frac{1}{1} + \frac{1}{3^r} &+ \ldots + \frac{1}{(2q-1)^r}\Big) - \Big( \frac{1}{(2q+1)^r} + \ldots + \frac{1}{(4q-1)^r}\Big) + \Big( \frac{1}{(4q+1)^r} + \ldots + \frac{1}{(6q-1)^r}\Big) \\
&- \Big( \frac{1}{(6q+1)^r} + \ldots + \frac{1}{(8q-1)^r}\Big) + \ldots \\
& \qquad \qquad = \frac{1}{(r-1)!q^r} \sum_{k=0}^{q-1}g^{(r-1)}(\frac{q-1-2k}{q}).
\end{split}
\end{equation}

As a representative sample of the sums that are obtained for various choices of $q$ and $r$, we have:

\begin{equation} \label{}
\begin{gathered}
\frac{1}{1} + \frac{1}{3} - \frac{1}{5} - \frac{1}{7} + \frac{1}{9} + \frac{1}{11} - \frac{1}{13} - \frac{1}{15} + \ldots = \frac{\pi \sqrt{2}}{4},\\
\frac{1}{1} + \frac{1}{3} + \frac{1}{5} - \frac{1}{7} - \frac{1}{9} - \frac{1}{11} + \frac{1}{13} + \frac{1}{15} +  \frac{1}{17} - \ldots = \frac{5 \pi}{12},\\
\frac{1}{1} + \frac{1}{3} + \frac{1}{5} + \frac{1}{7} - \frac{1}{9} - \frac{1}{11} - \frac{1}{13} - \frac{1}{15} +  \ldots = \pi \sqrt{2+ \sqrt{2}},\\
\frac{1}{1^3} + \frac{1}{3^3} - \frac{1}{5^3} - \frac{1}{7^3} + \frac{1}{9^3} + \frac{1}{11^3} - \frac{1}{13^3} - \frac{1}{15^3} + \ldots = \frac{3 \pi^3 \sqrt{2}}{128},\\
\frac{1}{1^3} + \frac{1}{3^3} + \frac{1}{5^3} - \frac{1}{7^3} - \frac{1}{9^3} - \frac{1}{11^3} + \frac{1}{13^3} + \frac{1}{15^3} +  \frac{1}{17^3} - \ldots = \frac{29 \pi^3}{864},\\
\frac{1}{1^5} + \frac{1}{3^5} - \frac{1}{5^5} - \frac{1}{7^5} + \frac{1}{9^5} + \frac{1}{11^5} - \frac{1}{13^5} - \frac{1}{15^5} + \ldots = \frac{\pi^5 57 \sqrt{2}}{24576},\\
\frac{1}{1^5} + \frac{1}{3^5} + \frac{1}{5^5} - \frac{1}{7^5} - \frac{1}{9^5} - \frac{1}{11^5} + \frac{1}{13^5} + \frac{1}{15^5} +  \frac{1}{17^5} - \ldots = \frac{1225 \pi^5}{373248}.
\end{gathered}
\end{equation}

If we desire the analogous sums which include the even integers, we may argue as follows. Let $q \geq 2$ be a positive integer, and consider $I(\frac{q-2}{q}) + I(\frac{q-4}{q}) + \ldots + I(\frac{-q+4}{q}) + I(\frac{-q+2}{q})$. Ordering the terms as before, it may be checked that the left side will be $q^r$ times the sum of the $r$-th powers of the reciprocals of the even integers other than the multiples of $2q$. We therefore obtain

\begin{equation} \label{armdark2}
\begin{split}
\Big( \frac{1}{2^r} + \frac{1}{4^r} &+ \ldots + \frac{1}{(2q-2)^r}\Big) - \Big( \frac{1}{(2q+2)^r} + \ldots + \frac{1}{(4q-2)^r}\Big) + \Big( \frac{1}{(4q+2)^r} + \ldots + \frac{1}{(6q-2)^r}\Big) \\
&- \Big( \frac{1}{(6q+2)^r} + \ldots + \frac{1}{(8q-2)^r}\Big) + \ldots \\
& \qquad \qquad = \frac{1}{(r-1)!q^r} \sum_{k=1}^{q-1} g^{(r-1)}(\frac{q-2k}{q}).
\end{split}
\end{equation}

In order to include the multiples of $2q$, add $\frac{1}{(2q)^r}$ times the series $\Delta_r :=\frac{1}{1^r}-\frac{1}{2^r} + \frac{1}{3^r} - \frac{1}{4^r} + \ldots$. Multiplying both sides by $2^r$, we obtain

\begin{equation} \label{embeth}
\begin{split}
\Big( \frac{1}{1^r} + \frac{1}{2^r} &+ \ldots + \frac{1}{q^r}\Big) - \Big( \frac{1}{(q+1)^r} + \ldots + \frac{1}{(2q)^r}\Big) + \Big( \frac{1}{(2q+1)^r} + \ldots + \frac{1}{(3q)^r}\Big) \\
&- \Big( \frac{1}{(3q+1)^r} + \ldots + \frac{1}{(4q)^r}\Big) + \ldots \\
& \qquad \qquad = \frac{\Delta_r}{q^r} + \frac{2^r}{(r-1)!q^r} \sum_{k=1}^{q-1} g^{(r-1)}(\frac{q-2k}{q}).
\end{split}
\end{equation}

The formula holds for $q=1$ as well, provided that the final sum on the right side is taken to be empty and therefore 0. Note that $\Delta_1 = \ln 2$, but that no other closed-form values of $\Delta_r$ for $r$ odd are known (nor are any likely to be known soon, as they can be expressed in terms of the Riemann zeta function evaluated at $r$). For $r=1, q=2,3,4$ we obtain

\begin{equation} \label{}
\begin{gathered}
\frac{1}{1} + \frac{1}{2} - \frac{1}{3} - \frac{1}{4} + \frac{1}{5} + \frac{1}{6} - \frac{1}{7} - \frac{1}{8} + \ldots = \frac{\pi}{4} + \frac{\ln 2}{2},\\
\frac{1}{1} + \frac{1}{2} + \frac{1}{3} - \frac{1}{4} - \frac{1}{5} - \frac{1}{6} + \frac{1}{7} + \frac{1}{8} + \frac{1}{9} - \ldots = \frac{2\pi}{3\sqrt{3}} + \frac{\ln 2}{3},\\
\frac{1}{1} + \frac{1}{2} + \frac{1}{3} + \frac{1}{4} - \frac{1}{5} - \frac{1}{6} - \frac{1}{7} - \frac{1}{8} + \ldots = \frac{\pi(1+2\sqrt{2})}{8} + \frac{\ln 2}{4}.
\end{gathered}
\end{equation}

\vski

\ccases{halfstrip} (half strip) Let $W$ now be the semi-infinite strip $\{-1<Re(z)<1, Im(z) > 0\}$. If $x \in (-1,1)$ and $a = \aa + \bb i \in W$, we can calculate $\rho_{T_W}^{\aa+\bb i}(x)ds$ by using the map $f(z) = \sin (\frac{\pi}{2}z)$, which maps $W$ conformally onto $\HH = \{Im(z) > 0\}$, fixing $-1$ and $1$. A simple calculation shows $\sin (\frac{\pi}{2}(\aa + \bb i)= \cosh(\frac{\pi \bb}{2})\sin(\frac{\pi \aa}{2}) + i \sinh(\frac{\pi \bb}{2})\cos(\frac{\pi \aa}{2}),$ and Theorem \ref{massage} and the calculation in Example \ref{halfplane} combine to give

\begin{equation} \label{michelle}
\begin{split}
\rho_{T_W}^{\aa+\bb i}(x)ds & = \frac{1}{2} \Big(\frac{\sinh(\frac{\pi \bb}{2})\cos(\frac{\pi \aa}{2})\cos(\frac{\pi x}{2})}{\sinh^2(\frac{\pi \bb}{2})\cos^2(\frac{\pi \aa}{2}) + (\cosh(\frac{\pi \bb}{2})\sin(\frac{\pi \aa}{2}) - \sin (\frac{\pi x}{2}))^2} \Big)ds \\
& = \frac{1}{2} \Big(\frac{\sinh(\frac{\pi \bb}{2})\cos(\frac{\pi \aa}{2})\cos(\frac{\pi x}{2})}{\sinh^2(\frac{\pi \bb}{2}) + \sin^2(\frac{\pi \aa}{2}) + \sin^2(\frac{\pi x}{2}) -2\cosh(\frac{\pi \bb}{2}) \sin(\frac{\pi \aa}{2}) \sin (\frac{\pi x}{2})} \Big)ds,
\end{split}
\end{equation}

where the identity $\sinh^2(\frac{\pi \bb}{2})\cos^2(\frac{\pi \aa}{2}) + \cosh^2(\frac{\pi \bb}{2})\sin^2(\frac{\pi \aa}{2}) = \sinh^2(\frac{\pi \bb}{2}) + \sin^2(\frac{\pi \aa}{2})$ was used to simplify the denominator. On the other hand, we can calculate this density using reflection as well. We claim that

\begin{equation} \label{}
\begin{split}
\rho_{T_W}^{\aa+\bb i}(x)ds &= \rho_{T_\HH}^{\aa+\bb i}(x)ds - \Big(\rho_{T_\HH}^{\aa+\bb i}(2-x)ds +\rho_{T_\HH}^{\aa+\bb i}(-2-x)ds\Big) \\
& \qquad \qquad + \Big(\rho_{T_\HH}^{\aa+\bb i}(4+x)ds +\rho_{T_\HH}^{\aa+\bb i}(-4+x)ds\Big) - \ldots \\
& = \frac{ds}{\pi} \Big(\frac{\bb}{\bb^2 + (\aa - x)^2} - \Big(\frac{\bb}{\bb^2 + (\aa - (2-x))^2} + \frac{\bb}{\bb^2 + (\aa - (-2-x))^2}\Big) \\
& \qquad \qquad+ \Big(\frac{\bb}{\bb^2 + (\aa - (4+x))^2} + \frac{\bb}{\bb^2 + (\aa - (-4+x))^2}\Big) - \ldots \Big)
\end{split}
\end{equation}

In order to justify this, we note first that $T_W = \min (T_\HH, \tau(-1), \tau(1))$, where $\tau$ was defined in the previous example. If a Brownian path exits $W$ at $x$, then it also exits $\HH$ at $x$, but must also not have hit $\{Re(z) = \pm 1\}$ before leaving $W$. The likelihood of exiting $\HH$ at $x$ after first striking $\{Re(z) = 1\}$ is the same, by reflection, as the likelihood of exiting $\HH$ at the reflection of $x$ over $\{Re(z) = 1\}$, which is $2-x$. Thus we must subtract $\rho_{T_\HH}^{\aa+\bb i}(2-x)ds$, and by the symmetric argument must also subtract the quantity obtained by reflection over $\{Re(z) = -1\}$, which is $\rho_{T_\HH}^{\aa+\bb i}(-2-x)ds$. However, we have now twice subtracted the contribution from paths which hit both of $\{Re(z) = \pm 1\}$ before hitting $\{Im(z)=0\}$, and we must therefore add $\rho_{T_\HH}^{\aa+\bb i}(4+x)ds +\rho_{T_\HH}^{\aa+\bb i}(-4+x)ds$; however, we have again overcompensated, and therefore must reflect again and subtract, and so forth (this argument can be made rigorous by adapting the methods given in Example \ref{strip}). Equating the two values obtained for the density, we obtain the identity

\begin{equation} \label{}
\begin{split}
\frac{1}{\pi} &\Big(\frac{\bb}{\bb^2 + (\aa - x)^2} - \Big(\frac{\bb}{\bb^2 + ((\aa+x) - 2)^2} + \frac{\bb}{\bb^2 + ((\aa+x) + 2)^2}\Big) \\
& \qquad + \Big(\frac{\bb}{\bb^2 + ((\aa - x) - 4)^2} + \frac{\bb}{\bb^2 + ((\aa-x) + 4)^2}\Big) - \ldots \Big) \\
& \qquad \qquad = \frac{1}{2} \Big(\frac{\sinh(\frac{\pi \bb}{2})\cos(\frac{\pi \aa}{2})\cos(\frac{\pi x}{2})}{\sinh^2(\frac{\pi \bb}{2}) + \sin^2(\frac{\pi \aa}{2}) + \sin^2(\frac{\pi x}{2}) -2\cosh(\frac{\pi \bb}{2}) \sin(\frac{\pi \aa}{2}) \sin (\frac{\pi x}{2})} \Big),
\end{split}
\end{equation}

where $x, \aa \in (-1,1)$, and $\bb > 0$. This identity is symmetric in $\aa$ and $x$, and if we set for instance $x=0$ we get

\begin{equation} \label{rach}
\begin{split}
\frac{1}{\pi} &\Big(\frac{\bb}{\bb^2 + \aa^2} - \Big(\frac{\bb}{\bb^2 + (\aa - 2)^2} + \frac{\bb}{\bb^2 + (\aa+ 2)^2}\Big) + \Big(\frac{\bb}{\bb^2 + (\aa - 4)^2} + \frac{\bb}{\bb^2 + (\aa + 4)^2}\Big) - \ldots \Big) \\
& \qquad \qquad = \frac{1}{2} \Big(\frac{\sinh(\frac{\pi \bb}{2})\cos(\frac{\pi \aa}{2})}{\sinh^2(\frac{\pi \bb}{2}) + \sin^2(\frac{\pi \aa}{2})} \Big).
\end{split}
\end{equation}

Setting $\aa=0$ as well gives

\begin{equation} \label{}
\frac{1}{\pi} \Big(\frac{1}{\bb} - \frac{2\bb}{\bb^2 + 2^2} + \frac{2\bb}{\bb^2 + 4^2} - \frac{2\bb}{\bb^2 + 6^2} + \frac{2\bb}{\bb^2 + 8^2} - \ldots \Big) = \frac{1}{2\sinh(\frac{\pi \bb}{2})}.
\end{equation}

If we subtract $\frac{1}{\pi \bb}$ from both sides, divide both sides by $\bb$, take the limit as $\bb \lar 0$, and simplify, we obtain

\begin{equation} \label{}
-1 + \frac{1}{2^2} - \frac{1}{3^2} + \frac{1}{4^2} - \ldots = \frac{-\pi^2}{12},
\end{equation}

which is easily seen to be equivalent to Euler's Basel sum \rrr{jill}. Returning to \rrr{rach}, if we assume $\aa \neq 0$, divide both sides by $\bb$, and take the limit as $\bb \lar 0$, we get

\begin{equation} \label{clea}
\frac{1}{\aa^2} - \Big(\frac{1}{(2-\aa)^2} + \frac{1}{(2+\aa)^2} \Big) + \Big(\frac{1}{(4-\aa)^2} + \frac{1}{(4+\aa)^2} \Big) - \ldots = \frac{\pi^2\cos(\frac{\pi \aa}{2})}{4 \sin^2(\frac{\pi \aa}{2})}.
\end{equation}



As before, the identity \rrr{clea} may be differentiated (and integrated) to obtain new identities if desired, and furthermore similar sums to \rrr{armdark} can be deduced. For example, taking an even number of derivatives and setting $x=0$ gives a sum equivalent to $\sum_{n=1}^{\ff} \frac{(-1)^n}{n^{2m-1}}$ for integer $m$, while taking an odd number of derivatives and letting $x \lar 1$ reduces to a sum equivalent to $\sum_{n=1}^\ff \frac{1}{n^{2m}}$ for integer $m$.

\vski

Let us now calculate the exit distribution from $W$ on the rays $\{Re(z) = \pm 1, Im(z) > 0\}$. The identity \rrr{michelle} holds upon replacing $x$ with $\pm 1 + yi$, provided that we recall that the $\cos(\frac{\pi x}{2})$ in the numerator came from the $|f'(z)|$ term in Theorem \ref{massage}, and therefore must be replaced by $|\cos(\frac{\pi}{2}(\pm 1 + yi))|$. We may also use $|\cos(\frac{\pi}{2}(\pm 1 + yi))| = \sinh (\frac{\pi}{2}y)$ for $y>0$ and $\sin(\frac{\pi}{2}(\pm 1 + yi)) = \pm \cosh (\frac{\pi}{2}y)$, and we obtain

\begin{equation} \label{mich2}
\rho_{T_W}^{\aa+\bb i}(\pm 1 + yi)ds  = \frac{1}{2} \Big(\frac{\sinh(\frac{\pi \bb}{2})\cos(\frac{\pi \aa}{2})\sinh(\frac{\pi y}{2})}{\sinh^2(\frac{\pi \bb}{2}) + \sin^2(\frac{\pi \aa}{2}) + \cosh^2(\frac{\pi y}{2}) \mp2\cosh(\frac{\pi \bb}{2}) \sin(\frac{\pi \aa}{2}) \cosh (\frac{\pi y}{2})} \Big)ds.
\end{equation}

This density can also be obtained by reflection, in several different ways. For example, if $W' = \{-1 < Re(z) < 1\}$ is the infinite strip of Example \ref{strip}, then $\rho_{T_W}^{\aa+\bb i}(\pm 1 + yi)ds = \rho_{T_{W'}}^{\aa+\bb i}(\pm 1 + yi)ds - \rho_{T_{W'}}^{\aa+\bb i}(\pm 1 - yi)ds$. This is because a Brownian path that leaves $W'$ at $\pm 1 + yi$ will also leave $W$ at that point, provided it does not first strike the real axis; however, the paths that strike the real axis and then proceed to $\pm 1 + yi$ will contribute the same probability, by reflection, as those that leave $W'$ at $\pm 1 - yi$. On the other hand, if $\tilde \tau(b) = \inf\{t \geq 0: Re(B_t) = b \mbox{ or } Im(B_t)=0\}$, then arguing similarly as in Example \ref{strip} we have

\begin{equation} \label{tats22}
\begin{gathered}
\rho^{\aa+\bb i}_{T_W}(1+yi)ds = \rho^{\aa+\bb i}_{\tilde \tau(1)}(1+yi)ds - \rho^{\aa+\bb i}_{\tilde \tau(-3)}(-3+yi)ds + \rho^{\aa+\bb i}_{\tilde \tau(5)}(5+yi)ds - \rho^{\aa+\bb i}_{\tilde \tau(-7)}(-7+yi)ds + \ldots, \\
\rho^{\aa+\bb i}_{T_W}(-1+yi)ds = \rho^{\aa+\bb i}_{\tilde \tau(-1)}(1+yi)ds - \rho^{\aa+\bb i}_{\tilde \tau(3)}(3+yi)ds + \rho^{\aa+\bb i}_{\tilde \tau(-5)}(-5+yi)ds - \rho^{\aa+\bb i}_{\tilde \tau(7)}(7+yi)ds + \ldots.
\end{gathered}
\end{equation}

It might seem as though here we have fertile ground for other identities, however it is easy to see, again by reflection, that $\rho^{\aa+\bb i}_{\tilde \tau(b)}(b+yi)ds = \rho^{\aa+\bb i}_{\tau(b)}(b+yi)ds - \rho^{\aa+\bb i}_{\tau(b)}(b-yi)ds$, where $\tau$ is the stopping time defined in Example \ref{strip}. Thus, any identities obtained here could just be obtained directly from \rrr{nosering}, evaluated at the proper values for $y$.

\vski

\ccases{Cminusslit} ($\CC \backslash [-1,1]$) Consider again the map $f(z) = \sin (\frac{\pi}{2}z)$ from the previous example, but this time let us examine the projection of the stopping time $T_\HH$ under this map. As stated before, this function maps the domain $\{-1<Re(z)<1, Im(z) > 0\}$ conformally onto $\HH$, and it does so by taking $\{Re(z)=-1,Im(z)>0\}$ onto $(-\ff,-1)$, $[-1,1]$ onto itself, and $\{Re(z)=1,Im(z)>0\}$ onto $(1,\ff)$; this can be verified by noting that $\sin(\frac{\pi}{2}(\pm 1+yi)) = \pm \cosh (\frac{\pi}{2} y)$. Schwarz reflection tells us then that $f$ maps $\{-3<Re(z)<-1, Im(z) > 0\}$ and $\{1<Re(z)<3, Im(z) > 0\}$ conformally onto $\{Im(z)<0\}$, and then $\{-5<Re(z)<-3, Im(z) > 0\}$ and $\{3<Re(z)<5, Im(z) > 0\}$ again conformally onto $\{Im(z)>0\}$, with every point on $\RR$ of course being mapped to a point in $[-1,1]$. It follows that the projection of $T_\HH$ under this map will be the first hitting time of the set $[-1,1]$, and we will abuse notation somewhat to refer to this stopping time as $T_{[-1,1]}$. If $x \in [-1,1]$, then the point $\sin(\frac{\pi}{2}x)$ will have preimages at $x, \pm 2 - x, \pm 4 + x, \pm 6 - x, \ldots$. Applying Theorem \ref{massage}, it follows that, for $\aa+\bb i \in \HH$, we have

\begin{equation} \label{yi}
\begin{split}
(\frac{\pi}{2} \cos(\frac{\pi}{2}x)) \rho^{\sin(\frac{\pi}{2}(\aa +\bb i))}_{T_{[-1,1]}} (\sin(\frac{\pi}{2}x)) ds & = \rho^{\aa + \bb i}_{T_\HH}(x)ds + \rho^{\aa + \bb i}_{T_\HH}(2-x)ds + \rho^{\aa + \bb i}_{T_\HH}(-2-x)ds\\
& \qquad + \rho^{\aa + \bb i}_{T_\HH}(4+x)ds + \rho^{\aa + \bb i}_{T_\HH}(-4+x)ds + \ldots \\
& = \frac{ds}{\pi}\Big(\frac{\bb}{\bb^2 + (\aa-x)^2} + \frac{\bb}{\bb^2 + (\aa-(2-x))^2} + \frac{\bb}{\bb^2 + (\aa-(-2-x))^2}\\
& \qquad + \frac{\bb}{\bb^2 + (\aa-(4+x))^2} + \frac{\bb}{\bb^2 + (\aa-(-4+x))^2} + \ldots \Big)
\end{split}
\end{equation}

Note that the term $(\frac{\pi}{2} \cos(\frac{\pi}{2}x))$ comes from the $|f'|$ in Theorem \ref{massage}, which by periodicity is equal at all preimages of $\sin(\frac{\pi}{2}x)$. The identity \rrr{yi} is similar in spirit to the identity obtained in Example \ref{punctdisk}, as every term on the right side of \rrr{yi} corresponds to a different homotopy class of Brownian curves hitting $[-1,1]$, some from above and some from below. In order to obtain the value for the right side we need to calculate $\rho_{T_{[-1,1]}}$ in a different way, and we can argue as follows. The M\"obius transformation $\phi(z) = \frac{1+z}{1-z}$ maps $[-1,1]$ to $[0,+\ff]$, and therefore maps $\hat \CC \backslash [-1,1]$ conformally onto $\hat \CC \backslash [0,+\ff]$, where $\hat \CC$ denotes the Riemann sphere. We can therefore use Theorem \ref{massage} to project the density for the first hitting time of $[-1,1]$ to the hitting time of $[0,+\ff)$. We obtain, for $\bar x = \sin(\frac{\pi}{2}x) \in [-1,1]$,

\begin{equation} \label{}
\rho^{\omega}_{T_{[-1,1]}}(\bar x) ds = \frac{2}{(1-\bar x)^2}\rho^{\phi(\omega)}_{T_{[0,+\ff)}}(\phi(\bar x)) ds.
\end{equation}

We can calculate $\rho_{T_{[0,+\ff)}}$ by projecting $\rho_{T_\HH}$ via the map $z \lar z^2$, which maps $\HH$ conformally onto $\CC \backslash [0,+\ff)$. The point $\phi(\bar x)$ will have two preimages, at $\pm \sqrt{\phi(\bar x)}$, and we obtain

\begin{equation} \label{yima}
\begin{split}
\rho^{\phi(\omega)}_{T_{[0,+\ff)}}(\phi(\bar x)) ds & = \frac{ds}{2\pi \sqrt{\phi(\bar x)}}\Big( \frac{Im \sqrt{\phi(\omega)}}{(Im \sqrt{\phi(\omega)})^2 + (Re \sqrt{\phi(\omega)}-\sqrt{\phi(\bar x)})^2} \\
& \qquad + \frac{Im \sqrt{\phi(\omega)}}{(Im \sqrt{\phi(\omega)})^2 + (Re \sqrt{\phi(\omega)}+\sqrt{\phi(\bar x)})^2} \Big),
\end{split}
\end{equation}

where the branch of the square root is chosen that takes values in $\HH \cup [0,+\ff)$. Combining \rrr{yi}-\rrr{yima} yields an identity, but this identity is fairly complex for arbitrary choice of $\aa + \bb i$. We can simplify considerably by taking $\aa$ to be an integer, and by periodicity we need only really consider $\aa=-1$ and $\aa = 0$. Let us begin with $\aa=-1$. As before we note that $\sin(\frac{\pi}{2}(-1+\bb i)) = -\cosh(\frac{\pi}{2} \bb)$, and thus $\phi(\sin(\frac{\pi}{2}(-1+\bb i))) = - \frac{\cosh(\frac{\pi}{2} \bb) - 1}{\cosh(\frac{\pi}{2} \bb) + 1}$, so that $\sqrt{\phi(\sin(\frac{\pi}{2}(-1+\bb i)))} = i \sqrt{\frac{\cosh(\frac{\pi}{2} \bb) - 1}{\cosh(\frac{\pi}{2} \bb) + 1}}$, a purely imaginary number. \rrr{yi} therefore simplifies to become

\begin{equation} \label{ima1}
\begin{split}
&\frac{\cos(\frac{\pi}{2}x)}{(1-\sin(\frac{\pi}{2}x))^2\sqrt{\frac{1+\sin(\frac{\pi}{2}x)}{1-\sin(\frac{\pi}{2}x)}}} \Big(\frac{\sqrt{\frac{\cosh(\frac{\pi}{2} \bb) - 1}{\cosh(\frac{\pi}{2} \bb) + 1}}}{\frac{\cosh(\frac{\pi}{2} \bb) - 1}{\cosh(\frac{\pi}{2} \bb) + 1} + \frac{1+\sin(\frac{\pi}{2}x)}{1-\sin(\frac{\pi}{2}x)}}\Big)ds\\
& \qquad \qquad \qquad = \frac{2}{\pi} \Big(\frac{\bb}{\bb^2 + (1+x)^2}ds + \frac{\bb}{\bb^2 + (3-x)^2}ds + \frac{\bb}{\bb^2 + (5+x)^2}ds + \frac{\bb}{\bb^2 + (7-x)^2}\Big)ds;
\end{split}
\end{equation}
\vski

note that the set of preimages in the sum in \rrr{yi} are symmetric around $-1$, thus the multiplicative factor of $2$ and the one-sided sum on the right side of \rrr{ima1}. Dividing both sides by $\bb$ and letting $\bb \searrow 0$, using $\lim_{\bb \searrow 0} \frac{\sqrt{\cosh(\frac{\pi}{2} \bb) - 1}}{\bb} = \frac{\pi}{2\sqrt{2}}$, gives

\begin{equation} \label{ima2}
\begin{split}
\frac{1}{(1+x)^2} + \frac{1}{(3-x)^2} + \frac{1}{(5+x)^2} + \frac{1}{(7-x)^2} + \ldots & = \frac{\pi^2 \cos(\frac{\pi}{2}x)}{8(1-\sin(\frac{\pi}{2}x))^2\Big(\frac{1+\sin(\frac{\pi}{2}x)}{1-\sin(\frac{\pi}{2}x)}\Big)^{3/2}} \\
& = \frac{\pi^2}{8(1+\sin(\frac{\pi}{2}x))}.
\end{split}
\end{equation}
If we take for instance $x=0$, we obtain

\begin{equation} \label{}
\frac{1}{1^2}+\frac{1}{3^2} + \frac{1}{5^2} + \frac{1}{7^2} + \ldots = \frac{\pi^2}{8},
\end{equation}

which is easily seen to be equivalent to Euler's Basel sum \rrr{jill}. In fact, it is not hard to see that \rrr{ima2} and \rrr{mei} are equivalent, with \rrr{ima2} simply being the derivative of \rrr{mei} (with $x = -a$).

\vski

Returning to \rrr{yi}, let us now see what happens when we take $\aa=0$. We need to calculate $\rho^{\sin(\frac{\pi}{2}(\bb i))}_{T_{[-1,1]}} (\sin(\frac{\pi}{2}x)) ds$, and we will use the same maps as before. We have $\sin(\frac{\pi}{2}\bb i) = i \sinh(\frac{\pi}{2}\bb i)$, and it may be checked that

\begin{equation} \label{}
\sqrt{\phi(i\sinh(\frac{\pi}{2}\bb))} = \sqrt{\frac{1+i\sinh(\frac{\pi}{2}\bb)}{1-i\sinh(\frac{\pi}{2}\bb)}} = \frac{1+i\sinh(\frac{\pi}{2}\bb)}{\sqrt{1+\sinh^2(\frac{\pi}{2}\bb)}} = \sech(\frac{\pi}{2}\bb) + i \tanh(\frac{\pi}{2}\bb).
\end{equation}

It may also be checked by simple manipulations that $\sqrt{\phi(\sin(\frac{\pi}{2}x))} = \sec(\frac{\pi}{2}x) + \tan(\frac{\pi}{2}x)$. These calculations, together with \rrr{yi}-\rrr{yima}, give the identity

\begin{equation} \label{}
\begin{split}
&\frac{\cos(\frac{\pi}{2}x)}{2(\sec(\frac{\pi}{2}x) + \tan(\frac{\pi}{2}x))(1-\sin(\frac{\pi}{2}x))^2} \Big( \frac{\tanh(\frac{\pi}{2}\bb)}{\tanh^2(\frac{\pi}{2}\bb) + (\sech(\frac{\pi}{2}\bb)-(\sec(\frac{\pi}{2}x) + \tan(\frac{\pi}{2}x)))^2} \\
& \qquad + \frac{\tanh(\frac{\pi}{2}\bb)}{\tanh^2(\frac{\pi}{2}\bb) + (\sech(\frac{\pi}{2}\bb)+(\sec(\frac{\pi}{2}x) + \tan(\frac{\pi}{2}x)))^2} \Big) \\
& = \frac{1}{\pi}\Big(\frac{\bb}{\bb^2 + x^2} + \frac{\bb}{\bb^2 + (2-x)^2} + \frac{\bb}{\bb^2 + (2+x)^2} + \frac{\bb}{\bb^2 + (4-x)^2} + \frac{\bb}{\bb^2 + (4+x)^2} + \ldots \Big).
\end{split}
\end{equation}

Divide both sides by $\beta$ and let $\beta \lar 0$ to obtain

\begin{equation} \label{}
\begin{split}
&\frac{1}{x^2} + \frac{1}{(2-x)^2} + \frac{1}{(2+x)^2} + \frac{1}{(4-x)^2} + \frac{1}{(4+x)^2} + \ldots \\
& \qquad =\frac{\pi^2 \cos(\frac{\pi}{2}x)}{4(\sec(\frac{\pi}{2}x) + \tan(\frac{\pi}{2}x))(1-\sin(\frac{\pi}{2}x))^2} \Big( \frac{1}{(1-(\sec(\frac{\pi}{2}x) + \tan(\frac{\pi}{2}x)))^2} + \frac{1}{(1+(\sec(\frac{\pi}{2}x) + \tan(\frac{\pi}{2}x)))^2} \Big) \\
& \qquad = \frac{\pi^2 \cos(\frac{\pi}{2}x)(1+(\sec(\frac{\pi}{2}x) + \tan(\frac{\pi}{2}x))^2)}{8\tan^2(\frac{\pi}{2}x)(\sec(\frac{\pi}{2}x) + \tan(\frac{\pi}{2}x))^3(1-\sin(\frac{\pi}{2}x))^2} \\
& \qquad = \frac{\pi^2 \cos(\frac{\pi}{2}x)(1+(\sec(\frac{\pi}{2}x) + \tan(\frac{\pi}{2}x))^2)}{8\sin^2(\frac{\pi}{2}x)(\sec(\frac{\pi}{2}x) + \tan(\frac{\pi}{2}x))}
\end{split}
\end{equation}

As with the earlier examples, a number of other identities can be deduced by evaluating the sum at particular values of $x$ and by differentiation. For example, taking the limit as $x \lar 1$ reduces to $\sum_{n=1}^{\ff} \frac{1}{(2n-1)^2} = \frac{\pi^2}{8}$, and subtracting $\frac{1}{x^2}$ from both sides and letting $x \lar 0$ leads to $\sum_{n=1}^{\ff} \frac{1}{n^2} = \frac{\pi^2}{6}$. Similarly, taking the suitable number of derivatives and letting $x$ tend to either $1$ or $0$ allow one to calculate the values of $\sum_{n=1}^\ff \frac{1}{n^{2m}}$ for integer $m$.

\vski

\ccases{rectangle} (rectangle) Let $W = \{-1<Re(z)<1,-k<Im(z)<k\}$ for some $k>0$. Let us calculate $\rho^{\aa+\bb i}_{T_W}(1+yi)ds$ in two different ways by reflection. Let us first note that, arguing as in Example \ref{halfstrip}, if $W' = \{-1 < Re(z) < 1\}$ then

\begin{equation} \label{aunt13}
\begin{split}
\rho_{T_W}^{\aa+\bb i}(1+yi)ds &= \rho_{T_{W'}}^{\aa+\bb i}(1+yi)ds - \Big(\rho_{T_{W'}}^{\aa+\bb i}(1+(2k-y)i)ds +\rho_{T_{W'}}^{\aa+\bb i}(1+ (-2k-y)i)ds\Big) \\
& \qquad \qquad + \Big(\rho_{T_{W'}}^{\aa+\bb i}(1+(4k+y)i)ds +\rho_{T_{W'}}^{\aa+\bb i}(1+ (-4k+y)i)ds\Big) - \ldots
\end{split}
\end{equation}

On the other hand, if $\hat \tau(b) = \inf\{t \geq 0: Re(B_t) = b \mbox{ or } Im(B_t)= \pm k\}$, then arguing similarly as in Example \ref{strip} we have

\begin{equation} \label{tats3}
\rho^{\aa+\bb i}_{T_W}(1+yi)ds = \rho^{\aa+\bb i}_{\hat \tau(1)}(1+yi)ds - \rho^{\aa+\bb i}_{\hat \tau(-3)}(-3+yi)ds + \rho^{\aa+\bb i}_{\hat \tau(5)}(5+yi)ds - \rho^{\aa+\bb i}_{\hat \tau(-7)}(-7+yi)ds + \ldots
\end{equation}

The two series given in \rrr{aunt13} and \rrr{tats3} must therefore be equal, with the terms in \rrr{aunt13} given by $\rrr{is3}$ and the terms in \rrr{tats3} given by an appropriate scaling, rotation, and translation of the corresponding density \rrr{michelle} for the domain $W''=\{-1<Re(z)<1, Im(z)>0\}$ considered in Example \ref{halfstrip}, since $\hat \tau(b)$ is the exit time of the domain bounded by $\{Re(z) = b\}$ and $\{Im(z) = \pm k\}$ containing $\aa+\bb i$, and this domain is clearly conformally equivalent to $W''$. Suffice it to say that the resulting identity is quite complex and not particularly illuminating. Let us therefore simplify things by setting $\aa=\bb=0$; the identity which results is then

\begin{equation} \label{}
\begin{split}
\sech& (\frac{\pi}{2} y) -\Big(\sech (\frac{\pi}{2} (2k-y)) + \sech (\frac{\pi}{2} (-2k-y))\Big) + \Big(\sech (\frac{\pi}{2} (4k+y)) + \sech (\frac{\pi}{2} (-4k+y))\Big) - \ldots \\
& = \frac{2\cos(\frac{\pi y}{2k})}{k}\Big(\frac{\sinh(\frac{\pi}{2k})}{\sinh^2(\frac{\pi}{2k}) + \sin^2(\frac{\pi y}{2k})} - \frac{\sinh(\frac{3\pi}{2k})}{\sinh^2(\frac{3\pi}{2k}) + \sin^2(\frac{\pi y}{2k})} + \frac{\sinh(\frac{5\pi}{2k})}{\sinh^2(\frac{5\pi}{2k}) + \sin^2(\frac{\pi y}{2k})} - \ldots \Big)
\end{split}
\end{equation}

Taking $y=0$, using the fact that $\sech$ is an even function, gives the identity

$$
\frac{1}{2} - \sech (\pi k) + \sech (2\pi k) - \sech (3\pi k) + \sech (4\pi k) - \ldots = \frac{1}{k} \Big( \csch (\frac{\pi}{2k}) - \csch (\frac{3\pi}{2k}) + \csch (\frac{5\pi}{2k}) - \csch (\frac{7\pi}{2k})+ \ldots\Big)
$$

\ccases{annulus} (annulus) Now let us consider the annulus $A_r=\{e^{-r}<|z|<e^{r}\}$, where $r>0$ is real. We map the strip $W=\{-1<Re(z)<1\}$ in Example \ref{strip} to $A_r$ via the function $f(z)=e^{rz}$, noting that $|\frac{d}{dz} e^z| = re^{r}$ on $\{Re(z)=1\}$ and $|\frac{d}{dz} e^z| = r e^{-r}$ on $\{Re(z)=-1\}$. We'll assume $a \in (e^{-r},e^r)$ is real. Theorem \ref{massage}, using \rrr{is3}, then gives

\begin{equation} \label{skyfall}
\begin{split}
\rho^a_{T_{A_r}}(e^{\pm r}e^{i \theta})ds &= \frac{1}{re^{\pm r}}\sum_{k=-\ff}^{\ff} \rho^{\frac{\ln a}{r}}_{T_W}(\pm 1+i(\frac{\th + 2\pi k}{r}))ds \\
& = \frac{1}{re^{\pm r}}\sum_{k=-\ff}^{\ff} \frac{\sech (\frac{\pi}{2}(\frac{\th + 2\pi k}{r}))(1-\tan^2(\frac{\pi\ln a}{4r}))}{4((1\mp\tan (\frac{\pi\ln a}{4r})\frac{1 - \tanh^2(\frac{\pi}{4}(\frac{\th + 2\pi k}{r}))}{1 + \tanh^2(\frac{\pi}{4}(\frac{\th + 2\pi k}{r}))})^2+(\tan (\frac{\pi\ln a}{4r})\frac{2\tanh(\frac{\pi}{4}(\frac{\th + 2\pi k}{r})) }{1 + \tanh^2(\frac{\pi}{4}(\frac{\th + 2\pi k}{r}))})^2)}ds.
\end{split}
\end{equation}

As with the punctured disk example, each term in the sum in \rrr{skyfall} corresponds to a homotopy class of curves leaving the annulus at that point. The formula is far simpler if we take $a=1$, in which case we have

\begin{equation} \label{}
\rho^a_{T_{A_r}}(e^{\pm r}e^{i \theta})ds = \frac{ds}{4 re^{\pm r}} \sum_{k=-\ff}^{\ff} \sech (\frac{\pi}{2}(\frac{\th + 2\pi k}{r}))
\end{equation}

We now show how applying Dynkin's formula with this expression yields an identity. If $h$ is harmonic on $A_r$, continuous on $\bar A_r$, then we obtain

\begin{equation} \label{chess}
\begin{split}
h(1) & = E_1[h(B_{T_{A_r}})] \\
&= \int_{0}^{2\pi} h(e^{-r + i\th}) \Big( \frac{1}{4 re^{-r}}\sum_{n=-\ff}^{\ff}\sech (\frac{\pi}{2}(\frac{\th + 2\pi k}{r}))\Big) (e^{-r} d\th) \\
& \qquad + \int_{0}^{2\pi} h(e^{r + i\th}) \Big( \frac{1}{4 r e^r}\sum_{n=-\ff}^{\ff}\sech (\frac{\pi}{2}(\frac{\th + 2\pi k}{r})) \Big) (e^{r}d\th) \\
& = \frac{1}{4r} \int_{-\ff}^{\ff} (h(e^{-r + i\th}) + h(e^{r + i \th}))\sech (\frac{\pi \th}{2r})d\th\\
& = \frac{1}{4} \int_{-\ff}^{\ff} (h(e^{-r + ir \th}) + h(e^{r + i r \th}))\sech (\frac{\pi \th}{2})d\th;
\end{split}
\end{equation}

note that the identity $ds=(e^{\pm r}d\th)$ was used on the curves $\{|z|=e^{\pm r}\}$. If we put in $h(z) = z$ and rearrange we obtain the identity

\begin{equation} \label{sweet}
\int_{-\ff}^{\ff} e^{i r \th} \sech(\frac{\pi \th}{2}) d\th = 2 \sech (r).
\end{equation}

It should be noted that the same identity with $r < 0$ is obtained by setting $h(z) = \bar z$. We have therefore derived the Fourier transform of the $\sech$ function.

\vski

\noindent {\bf Remark:} It may be tempting to look for other harmonic functions on the annulus in order to derive new identities from \rrr{chess}, in particular we might hope that the argument given here provides a general method for evaluating Fourier transforms involving the function $\sech$. However, the reader should be aware that a search for further identities using this exit distribution will lead to nothing else substantial. This is because for any harmonic function $h$ on $A_r$ an analytic function $g$ and real constant $C$ can be found such that $h = Re(g(z)) + C\log|z|$ (see for example \cite[Ex. III.3.4]{gamcomp}). $\log |z|$ in \rrr{chess} yields a triviality, and analytic functions on annuli are Laurent series in $z$, so that if we write $h-C\log|z| = \frac{g + \bar g}{2}$ we see that $h-C\log|z|$ can be expressed as a Laurent series in $z$ and $\bar z$. Applying \rrr{chess} to $h$ would therefore not yield anything more than would a term by term application with $h(z) = z^q$ or $h(z) = \bar z^q$, each of which give simply \rrr{sweet}.

\vski

It is also interesting to see what happens if we try to adapt the reflection technique from Example \ref{strip} to this case. If the analogous argument to that in Example \ref{strip} applied, we would have

\begin{equation} \label{thai}
\rho^a_{T_{A_r}}(e^{r}e^{i \theta})ds = \rho^a_{\hat \tau(e^r)}(e^{r}e^{i \theta})ds -\rho^a_{\hat \tau(e^{-r},e^r)}(e^{r}e^{i \theta})ds + \rho^a_{\hat \tau(e^r,e^{-r},e^r)}(e^{r}e^{i \theta})ds - \rho^a_{\hat \tau(e^{-r},e^r,e^{-r},e^r)}(e^{r}e^{i \theta})ds + \ldots
\end{equation}

where $\hat \tau(m_1,m_2,\ldots, m_k)$ would be the first time the Brownian motion has hit all of the curves $\{|z|=m_1\}, \{|z|=m_2\}, \ldots, \{|z|=m_k\}$ in order. The quantities on the right side of \rrr{thai} are easy to find using reflection (note that reflection over $\{|z|=m\}$ is given by the function $z \lar \frac{m^2}{\bar z}$) and the calculations in Example \ref{disk}, and if we consider the simplest case, when $\th = 0$, $a=1$, we would obtain the identity

\begin{equation} \label{coco} \nonumber
\frac{1}{4re^{r}} \sum_{k=-\ff}^{\ff} \sech (\frac{\pi^2 k}{r}) = \frac{1}{2\pi e^r} \Big( \frac{e^r + 1}{e^r-1} - \frac{1+e^{-3r}}{1-e^{-3r}} + \frac{e^{5r} + 1}{e^{5r}-1} - \frac{1+e^{-7r}}{1-e^{-7r}} + \ldots \Big)
\end{equation}

However, this sum cannot be valid, since the sum on the right does not converge (the terms inside the parentheses approach 1). The reason that the argument fails in this case is that the densities $\rho^a_{\hat \tau(e^{-r},e^r, \ldots,e^r)}$ and $\rho^a_{\hat \tau(e^{r},e^{-r}, \ldots,e^r)}$ do not approach $0$ as the length of the sequence goes to infinity, but rather approach the uniform density on the circle. On the other hand, the analogous statement to Lemma \ref{mang} will still hold, so it might be interesting to see whether any sense can be made of \rrr{coco} or the reflection argument.

\vski

\ccases{wind1} We now consider a stopping time which is not the exit time of a domain. Start a Brownian motion at $1$, and let $\tau_r = \inf\{t: arg(B_t) = \pm r\pi\}$, with the branch of the argument chosen so that $arg(B_0)=0$ a.s. We can calculate the distribution of $B_{\tau_r}$, and arrive at

\begin{equation} \label{mae2}
\rho^1_{\tau_r}(ye^{\pm r\pi i})ds = \frac{ds}{2r\pi y(y^{\frac{1}{2r}}+y^{\frac{-1}{2r}})},
\end{equation}

when $r$ is not an integer. When $r$ is an integer, then $e^{r\pi i} = e^{-r\pi i}$, so doubling \rrr{mae2} we have

\begin{equation} \label{}
\rho^1_{\tau_r}(ye^{r\pi i})ds = \frac{ds}{r\pi y(y^{\frac{1}{2r}}+y^{\frac{-1}{2r}})}.
\end{equation}

To see this, we project the density found earlier of the half-plane via the map $f(z)=z^{2r}$ with $|f'(yi)| = 2ry^{2r-1}$ to get

\begin{equation} \label{}
\rho^1_{\tau_r}(y^{2r}e^{\pm \pi r i}) ds = \frac{ds}{\pi(1+y^2)2ry^{2r-1}}.
\end{equation}

Replace $y$ with $y^{\frac{1}{2r}}$ to arrive as claimed at

\begin{equation} \label{}
\rho^1_{\tau_r}(y e^{\pm \pi r i}) ds = \frac{ds}{2\pi r(1+y^{\frac{2}{2r}})y^{\frac{2r-1}{2r}}} = \frac{ds}{2\pi ry(y^{\frac{1}{2r}}+y^{\frac{-1}{2r}})}.
\end{equation}

This density admits the antiderivative $\frac{1}{\pi} \tan^{-1}(y^{\frac{1}{2r}})$, and this allows us in particular to note that

\begin{equation} \label{}
P_1(|B_{\tau_r}| \in (0,\eps)) = \frac{2}{\pi} \tan^{-1}(\eps^{\frac{1}{2r}}) = P_1(|B_{\tau_r}| \in (\frac{1}{\eps}, \ff)),
\end{equation}

and this probability approaches $\frac{1}{2}$ as $r$ approaches $\ff$. Thus, as $r \lar \ff$, the distribution of $|B_{\tau_r}|$ approaches $\frac{1}{2}(\dd_0 + \dd_\ff)$ in distribution. This may seem surprising at first, especially the mass accumulating at $0$; but in fact it is to be expected since any mass at infinity must correspond to a mass at 0 due to the fact that $\frac{1}{B_t}$ is a (time-changed) Brownian motion. Evidently when $B_t$ comes close to $0$ it winds many times around the origin, analogously to how a one-dimensional Brownian motion hits $0$ infinitely often in any neighborhood of a visit to $0$, and this results in a large change in argument and hence a high likelihood of being near 0 at the first time attaining a prescribed argument.

\vski

We also remark briefly that the symmetry imposed upon $\tau_r$ can easily be dispensed of. In other words, we may let $\tau_{r_1,r_2} = \inf\{t: arg(B_t) = r_1\pi \mbox{ or } -r_2\pi\}$ for $r_1, r_2 >0$, and calculate

\begin{equation} \label{lamchop}
\begin{gathered}
\rho^1_{\tau_{r_1,r_2}}(ye^{r_1 \pi i})ds = \frac{\cos \th \quad \! \! \! ds}{\pi (r_1+r_2)y^{1-\frac{1}{r_1+r_2}}(\cos ^2 \th+(y^{\frac{1}{r_1+r_2}} - \sin \th)^2)},\\
\rho^1_{\tau_{r_1,r_2}}(ye^{-r_2 \pi i})ds = \frac{\cos \th \quad \! \! \! ds}{\pi (r_1+r_2)y^{1-\frac{1}{r_1+r_2}}(\cos ^2 \th+(-y^{\frac{1}{r_1+r_2}} - \sin \th)^2)},
\end{gathered}
\end{equation}

where $\th = \frac{\pi}{2} \Big(\frac{r_2-r_1}{r_2+r_1}\Big)$. This is obtained by first noting by rotational invariance that we may obtain the distribution by starting a Brownian motion at $e^{\frac{\pi}{2}(r_2-r_1)i}$ and stopping it at the first time the argument reaches $\pm \frac{\pi}{2}(r_2+r_1)$, and then projecting the density from the right half-plane as before by the function $f(z)=z^{r_1+r_2}$. Note that if $e^{r_1 i} = e^{-r_2 i}$ then these two quantities must be added in order to find the density of $B_{\tau_{r_1,r_2}}$ on the ray $\{arg(z) = r_1\} = \{arg(z) = -r_2\}$.

\vski

We may also stop the Brownian motion at a prescribed argument. Let $\hat \tau_r = \inf\{t: arg(B_t) = r\}$, where we will assume $r>0$. If we let $W = \{Im(z) < r\}$, then the stopping time $T_W$ projects to $\hat \tau_r$ under the exponential map $z \lar e^z$. Theorem \ref{massage} gives

\begin{equation} \label{}
\rho_{\hat \tau_r}^1(ye^{ir})ds = \frac{1}{y}\rho^0_{T_W}(\ln y + ri) ds = \frac{1}{y\pi} \frac{r}{r^2 + (\ln y)^2} ds.
\end{equation}

It is straightforward to find the corresponding distribution function if desired, and to note that the accumulation of mass at $0$ and $\ff$ holds for this example as well; that is, $P_1(|B_{\hat \tau_r}| \in (0,\eps)) = P_1(|B_{\hat \tau_r}| \in (\frac{1}{\eps}, \ff)) \lar \frac{1}{2}$ as $r \lar \ff$. It may also be noted that reflection may be used in order to express $\rho_{\tau_r}$ or $\rho^1_{\tau_{r_1,r_2}}$ as an alternating infinite sum of terms of the form $\rho_{\hat \tau_r}^1$, exactly in the same manner as in Exercise \ref{strip}. The identities obtained in this manner are precisely the same as in Example \ref{strip}, since the exponential function will map the properly chosen strip or half-plane into the winding stopping times.

\vski

\ccases{what} Now let $\hat \tau = \inf_{t \geq 0} (B_t \in (-1,1), \{B_s\}_{0 \leq s \leq t} \cup [0,B_t] \mbox{ is not homotopic to a point} )$; that is, $t$ is the first time that $B_t$ lies on $(-1,1)$ simultaneously with the curve traced by $B_s$ up to time $t$ being wound around at least one of $-1$ and $1$. We will calculate the distribution of $B_{\hat \tau}$. First we set $\tau = \inf\{t: B_t \in (-\ff,-1] \cup [1,+\ff) \}$. We will show that $\rho^{0}_\tau(w) ds = \frac{1}{\pi |w| \sqrt{w^2-1}}$, and then project this density to $\hat \tau$.  Let $\tau_1 = \inf\{t: B_t \in (-\ff,+\ff) \}$ and $\tau_0 = \inf\{t: B_t \in [0,+\ff) \}$. Using the map $f(z)=z^2$ and Theorem \ref{massage}, we project the density for $\tau_1$ of Example \ref{halfplane} via the conformal map $f(z)=z^2$ to get $\rho^{(-1)}_{\tau_2}(v) ds = \frac{1}{\pi (1+v)\sqrt{v}}$ (the $-1$ is placed in parentheses to prevent any confusion with an inverse map). We can now project this density via the transformation $w=\phi(v)=\frac{1+v}{1-v}$ to obtain

\begin{equation} \label{}
\begin{split}
\rho^{(0)}_\tau(v) ds & = \rho^{(-1)}_{\tau_1}(\frac{w-1}{w+1}) (\phi^{-1})'(w)  ds \\
& = \frac{1}{\pi \sqrt{\frac{w-1}{w+1}}(1+\frac{w-1}{w+1})}\times \frac{2}{(w+1)^2} ds \\
& = \frac{1}{\pi |w|\sqrt{w^2-1}} ds,
\end{split}
\end{equation}

as claimed. We now will project the density of $\tau$ to $\hat \tau$ using the entire function $f(z) = \sin(\frac{\pi}{2}z)$. To see that this does the job, note that $\sin(\frac{\pi}{2}(-1+yi))=-\cosh(\frac{\pi}{2}y)$ and $\sin(\frac{\pi}{2}(1+yi))=\cosh(\frac{\pi}{2}y)$. $f$ therefore maps the boundary of the half-infinite strip $\{Im(z)>0, -1 < Re(z) <1\}$ injectively onto the boundary of $\{Im(z)> 0\}$. The argument principle allows us to conclude that $f$ maps $\{Im(z)>0, -1 < Re(z) <1\}$ conformally onto $\{Im(z)> 0\}$. The Schwarz reflection principle now informs us that $f$ maps $\{Im(z)>0, -3 < Re(z) <-1\}$ conformally onto $\{Im(z)< 0\}$, with $(-3,-1)$ being mapped to $(-1,1)$. Reflecting in this manner to the left and right, as well as below to $\{Im(z)<0, -1 < Re(z) <1\}$ and thence left and right to fill out the plane, the structure of $f$ can be understood; in particular we see that a closed curve $\CCC$ beginning and ending at $0$ (but otherwise not touching $(-1,1)$) is not homotopic to a point if and only if there is a curve $\ga$ traveling from $0$ to some non-zero even integer such that $f(\ga) = \CCC$. This tells us that $\tau$ is mapped to $\hat \tau$ under $f$. If $w \in (-1,1)$, then the preimages under $f$ of $w$ on $(-\ff,-1] \cup [1,+\ff)$ are all points of the form $(4n-2) - \frac{2}{\pi}\sin^{-1} w$ for integer $n$ or $\frac{2}{\pi}\sin^{-1} w+4n$ for integer $n \neq 0$. Note that the value of $|f'|$ at these points is $\frac{\pi}{2}\cos(\sin^{-1} w) = \frac{\pi}{2}\sqrt{1-w^2}$. Applying our theorem therefore gives

\begin{equation} \label{}
\begin{split}
\rho^0_{\hat \tau}(w)ds = &\frac{2ds}{\pi^2 \sqrt{1-w^2}} \Big( \sum_{\substack{n=-\ff \\ n\neq 0}}^\ff \frac{1}{|\frac{2}{\pi}\sin^{-1} w+4n|\sqrt{(\frac{2}{\pi}\sin^{-1} w+4n)^2-1}} \\
&+ \sum_{n=-\ff}^\ff \frac{1}{|\frac{2}{\pi}\sin^{-1} w-(4n+2)|\sqrt{(\frac{2}{\pi}\sin^{-1} w-(4n+2))^2-1}}\Big)\\
= & \frac{2ds}{\pi^2 \sqrt{1-w^2}} \sum_{\substack{n=-\ff \\ n\neq 0}}^\ff \frac{1}{|\frac{2}{\pi}\sin^{-1} w+2n|\sqrt{(\frac{2}{\pi}\sin^{-1} w+2n)^2-1}}.
\end{split}
\end{equation}

Somewhat similarly to earlier examples, each term in the sum here corresponds to homotopy classes of Brownian curves, with each term corresponding to two different classes: one in which $B_t$ approaches $B_\tau \in (-\ff,-1] \cup [1,+\ff)$ from above, and one in which it approaches from below. Note that

\begin{equation} \label{}
\frac{d}{dw} \Big(\frac{-1}{\pi} \cot^{-1} \sqrt{(\frac{2}{\pi}\sin^{-1} w+2n)^2 - 1}\Big) = \frac{2ds}{\pi^2 \sqrt{1-w^2}(\frac{2}{\pi}\sin^{-1} w+2n)\sqrt{(\frac{2}{\pi}\sin^{-1} w+2n)^2-1}}.
\end{equation}

The distribution function is therefore explicitly calculable, if desired.

\section{Acknowledgements} The author is grateful for support from Australian Research Council Grants DP0988483 and DE140101201.

\bibliographystyle{alpha}
\bibliography{CAbib}

\end{document}